# ENDS IN FREE MINIMAL SPANNING FORESTS[1]

By Ádám Timár

*Indiana University*

We show that for a transitive unimodular graph, the number of ends is the same for every tree of the free minimal spanning forest. This answers a question of Lyons, Peres and Schramm.

**1. Introduction.** Let $G = (V, E)$ be a graph and let $\{U(e)\}_{e \in E}$ be independent uniform $[0, 1]$ random labels assigned to its edges. Define the *free minimal spanning forest* [FMSF or $\mathsf{FMSF}(G)$] of $G$ as the set of edges $e$ where $U$ does not attain its maximum at $e$ for any cycle containing $e$. The *wired minimal spanning forest* [WMSF or $\mathsf{WMSF}(G)$] is defined analogously: this is the set of edges $e$ where $U$ does not attain its maximum at $e$ for any cycle or any bi-infinite simple path containing $e$. It is easy to see that the random graphs defined this way are indeed spanning forests. When there is no need to specify whether we refer to the free or the wired forest, we call it MSF.

Minimal spanning forests of finite graphs have been studied for a long time. In this case, MSF is the tree that has minimum sum of weights over its edges, which make MSFs essential for certain optimization problems. The interest towards minimal spanning forests in infinite graphs arose because of their close connections to percolation. See the paper of Lyons, Peres and Schramm [2] for the history and further references.

In what follows, it is assumed that the infinite graph $G$ considered is locally finite, transitive and unimodular (e.g., it may be a finitely generated Cayley graph).

In [2] various properties of MSFs are proved. Several open questions are listed there, one of which is the following. Given a transitive unimodular graph $G$, is it true that if $\mathsf{FMSF}(G) \neq \mathsf{WMSF}(G)$, then every tree of

Received March 2005; revised May 2005.

[1]Supported in part by NSF Grant DMS-02-31224 and Hungarian National Foundation for Scientific Research Grant TO34475.

FMSF($G$) has infinitely many ends? The condition FMSF($G$) $\neq$ WMSF($G$) implies that at least one tree of FMSF($G$) has infinitely many ends. From Proposition 3.6 of [2] we know that FMSF($G$) $\neq$ WMSF($G$) can hold only if $G$ is nonamenable. In this case every component of WMSF($G$) has one end. Hence our question can be rephrased: does every tree of FMSF($G$) have infinitely many ends provided that one of them has infinitely many ends? Or, equivalently, are the (infinite) components of the FMSF of a nonamenable graph indistinguishable by the number of their ends? The indistinguishability of infinite clusters is known when the percolation process is insertion tolerant [3], and the same question about MSFs is one of the challenging questions about MSFs.

In this paper we shall prove that the trees of the FMSF all have the same number of ends. Note that if an invariant tree has infinitely many ends, then there are $2^{\aleph_0}$ ends. (For a proof one can use Lemma 2.2.) Hence, our result says that if FMSF($G$) $\neq$ WMSF($G$), then every FMSF-tree has uncountably many ends.

For its central role, let us give the definition: an *end* of a tree $T$ is an equivalence class on the set of (simple) paths in $T$, where two paths are equivalent iff their symmetric difference is finite. Suppose that $P$ is a representing path for an end $\xi$ and $x$ is some vertex of $P$. Let $C$ be the component of $T \setminus \{x\}$ that contains all but finitely many vertices of $P$. If $x$ could be chosen so that any path in $C$ differs from $P$ in only finitely many edges (i.e., if $C$ has only one end), then we say that $\xi$ is an *isolated end*.

We shall use the *mass transport principle* (MTP) several times, by the following simple corollary of it.

PROPOSITION 1.1. *Let $G$ be a unimodular transitive graph, and let $S$ be a random set of disjoint infinite subsets of $V(G)$, where the distribution of $S$ is invariant under the automorphisms of $G$. Then one cannot assign a finite nonempty subset to each (or at least one) set in $S$ in a way that is equivariant with $S$.*

A description of the mass transport principle can be found, for example, in [2].

We shall not always mention "almost always" when this is the case. We shall use without mentioning that MSF is ergodic (as proved in [2]); hence we do not distinguish between invariant events of positive probability or of probability 1.

We shall prove the following.

THEOREM 1.2. *If $G$ is a transitive, unimodular graph and FMSF($G$) is not equal to WMSF($G$), then every tree of FMSF has infinitely many ends.*

REMARK 1.3. Notice that the claim here is equivalent to saying that no WMSF-tree can coincide with an FMSF-tree. The number of WMSF-trees contained in an FMSF-tree $T$ is 1 more than the number of edges in $T \setminus \mathsf{WMSF}$ (in case of infinitely many, they are equal). Since one cannot choose a finite subset of $T$ in an invariant way, $T$ has infinitely many ends if and only if infinitely many WMSF-trees are contained in it. By Theorem 3.12 of [2] the converse of Theorem 1.2 is also true, because any tree of the WMSF has one end.

Let $G_p$ denote the graph formed by the set of edges with labels $< p$ and their endpoints. Let $G_p^*$ stand for the union of infinite components of $G_p$. Finally, given a configuration $\kappa$ of edge labels, denote by $\kappa_p$ the set of edges whose label is smaller than $p$ in $\kappa$; $\mathsf{WMSF}(\kappa)$ and $\mathsf{FMSF}(\kappa)$ stand for the respective spanning forests for the configuration $\kappa$.

**2. Infinitely many ends of FMSF-trees.** In this section we prove Theorem 1.2.

Here is a brief sketch of the proof. If $\mathsf{FMSF} \neq \mathsf{WMSF}$, then $p_c < p_u$. Choose a $p$ strictly between $p_c$ and $p_u$. We shall use a kind of "weak" insertion tolerance for the FMSF, that we define in the rest of this paragraph. Suppose that $A$ is an event of positive probability and $e$ is an edge such that, on $A$, the endpoints of $e$ are in distinct connected components of $G_p$. Define $A'$ as the event arising from $A$ when we multiply the label of $e$ by $p$. Then $\mathbf{P}[A'] > 0$. Further, $\mathsf{FMSF} \cap G_p$ on $A'$ is the union of $e$ and the $\mathsf{FMSF} \cap G_p$ of the corresponding configuration in $A$. Hence, if we denote the $\mathsf{FMSF}(G) \cap G_p$-trees that contain the endpoints of $e$ on $A$ by $T$ and $T'$, then $\mathsf{FMSF}(G) \cap G_p$ on $A'$ contains a tree that is the union of $e$, $T$ and $T'$. Also, it is easy to check that $T$ and $T'$ are infinite, and they belong to different trees of $\mathsf{FMSF}(G)$.

So the edge $e$ was "inserted" in the FMSF and the new event still has positive probability. The definition that we outlined here is actually simpler (though basically the same) as the one that we shall need in the proof.

Now, if we suppose that there is a tree in FMSF with finitely many ends, then "insert" a path between this tree and another FMSF-tree, to get an FMSF-tree with at least two ends, some of which are isolated. Such a tree cannot occur with positive probability, by a result proved in [3], giving a contradiction.

From now on, we *confine ourselves to nonamenable graphs*. As we have already mentioned, the hypothesis of Theorem 1.2 fails for amenable graphs. For nonamenable graphs, we will use repeatedly the fact that WMSF-trees have one end, by Theorem 3.12 in [2].

The following is Proposition 3.6 in [2].

LEMMA 2.1. *Let $G$ be a unimodular transitive graph. Then the following are equivalent:*

  (i) $\mathsf{FMSF}(G)$ *is not equal to* $\mathsf{WMSF}(G)$.
 (ii) $p_c(G) < p_u(G)$.

Another tool we need is the following lemma, which was first stated in [3].

LEMMA 2.2.  *Let $F$ be a random forest in a unimodular transitive graph $G$ whose distribution is invariant under the automorphism group of $G$. Then no tree of $F$ with infinitely many ends can have an isolated end.*

PROOF.  Otherwise we could assign vertex $x$ to each maximal 1-ended component of $F \setminus \{x\}$, giving an MTP contradiction.  □

We say that the components $C_1$ and $C_2$ of a subgraph of a graph $G$ are *connected* by a path $P$ if $C_1 \cup C_2 \cup P$ is connected and one endpoint of $P$ is in $C_1$ and the other one is in $C_2$. An *inner vertex* of a path $P$ is a vertex different from the endpoints of $P$.

PROOF OF THEOREM 1.2.  We prove by contradiction. Suppose that there is a tree of $\mathsf{WMSF}(G)$ that is also a tree of $\mathsf{FMSF}(G)$. (This is the negation of our claim by Remark 1.3.) Call trees with this property *lonely*. As a corollary of Proposition 3.1 of [1] (also repeated in [2]), the MSF intersects some infinite cluster of $G_p$ in an infinite component whenever $p > p_c$.

Fix some $p \in (p_c, p_u)$. Such a $p$ exists by Lemma 2.1. There exists a finite path $P$ in $G$ such that with positive probability the following hold:

  (i) $P$ connects two components $K_1$ and $K_2$ of $G_p^*$.
 (ii) The endpoint of $P$ in $K_1$ is a vertex from a lonely tree $T$ (and of course $T \cap K_1$ infinite).
(iii) No edge of $T \cap K_1$ is incident to any *inner* vertex of $P$.

Denote by $D$ the set of edges incident to some *inner* vertex of $P$. Hence (iii) says that $T \cap K_1 \cap D$ is empty.

Such a choice indeed exists because the countable union of events satisfying the first two conditions for some finite $P$ is just the event of having a lonely tree and more than one component in $G_p^*$ (which has probability 1). The last condition is fulfilled if we choose a path $P$ of minimal length (and satisfying the other conditions).

Fix $P$ and let $E$ be the event that (i), (ii) and (iii) hold. Define $K_1$ and $K_2$ as in the criteria. For a configuration $\kappa$ in $E$ we can perform some of the following transformations:

(1) Change the label $U(e)$ of each edge $e \in D \setminus P$ to a new label $p + (1-p)U(e)$.

(2) Change the label $U(f)$ of each edge $f \in P$ to a new label $pU(f)$.

Let $\kappa''$ be the configuration we get after applying (1) on $\kappa$; let $\kappa'$ be the one we get after applying both transformations. The sets of every $\kappa''$ and $\kappa'$ arising by these couplings also have positive measures, since $\mathbf{P}[E] > 0$. Now, the MSF in $\kappa''$ differs from $\mathsf{MSF}(\kappa)$ only by finitely many edges, since the values on the edges differ only at finitely many places (Lemma 3.15 of [2]). Since we do not decrease labels going from $\kappa$ to $\kappa''$, the only thing that could change from $\mathsf{FMSF}(\kappa) \setminus D \;(\supset T)$ to $\mathsf{FMSF}(\kappa'')$ is that a few edges become part of the FMSF (but only those in $D \setminus P$ might "fall out"). However, almost surely none of these new edges can connect $T$ with some other tree of $\mathsf{FMSF}(\kappa_p)$. Otherwise a tree of $\mathsf{FMSF}(\kappa_p)$ with more than one end (i.e., infinitely many ends) would contain $T$ and an isolated end in it, contradicting Lemma 2.2. Hence $T$ is in a lonely tree in $\kappa''$ as well as in $\kappa$. (We mention that in $\kappa''$ the lonely tree containing $T$ may have finitely many edges not in $T$.)

From $\kappa''$ we get $\kappa'$ by changing the labels of the edges in $P$, according to (2). Every edge of $P$ becomes part of $\kappa'_p$. They are also in $\mathsf{FMSF}(\kappa')$, because any cycle through any edge $f \in P$ intersects $E(G) \setminus \kappa'_p$. (By the construction of $\kappa''$, as we have already pointed out, no inner vertex of $P$ is incident to any edge in $\kappa'_p \setminus P$.) Similarly, every edge of $\mathsf{FMSF} \cap \kappa''_p$ is also in $\mathsf{FMSF} \cap \kappa'_p$, because any cycle that contains an edge in $P$ (i.e., an edge whose label was changed) also intersects $E(G) \setminus \kappa'_p$.

Now we get to the final contradiction. Let $S$ be the tree in $\mathsf{WMSF}(\kappa'') \cap K_2$ adjacent to an endpoint of $P$. It is easy to check that $|S \cap K_2|$ is infinite. Since $P$ belongs to $\mathsf{FMSF}(\kappa')$, the edges of $T \cup P \cup S$ are in a tree $T'$ of $\mathsf{FMSF}(\kappa')$ [by the conclusions of the previous paragraph and that $T \cup S \subseteq \mathsf{FMSF}(\kappa'')$]. Hence $T'$ has more than one end. Now, $T'$ either has finitely many ends, in which case $T' \setminus \mathsf{WMSF}(\kappa')$ is finite, giving an MTP contradiction, or $T'$ has infinitely many ends, in which case it has an isolated end (provided by $T$), but that is impossible by Lemma 2.2.

This completes the proof. □

REMARK 2.3. It is possible that our method can be used to prove indistinguishability of $\mathsf{FMSF}$-trees when $p_c < p_u$. The "weak" insertion tolerance at $p$ $(p_c < p < p_u)$ makes the arguments in [3] applicable in this setting, provided that the existence of distinguishable $\mathsf{FMSF}$-trees in $G$ would imply that there are distinguishable trees in $\mathsf{FMSF}(G) \cap G_p$. Unfortunately we could not prove this implication. Further interesting properties, such as relentless merging, would follow too. By this we mean that for $p_c < p_1 < p_2 < p_u$ any tree in $\mathsf{MSF} \cap G_{p_2}$ contains infinitely many trees from $\mathsf{MSF} \cap G_{p_1}$.

**Acknowledgments.** I thank Russell Lyons and Gábor Pete for their comments on the manuscript.

DEPARTMENT OF MATHEMATICS
INDIANA UNIVERSITY
BLOOMINGTON, INDIANA 47405-5701
USA
E-MAIL: atimar@indiana.edu
URL: mypage.iu.edu/~atimar